\documentclass[10pt]{article}
\usepackage[dvips]{graphics}
\usepackage[cp866]{inputenc}
\input epsf
\usepackage{graphicx}

\begin{document}

\title{Sobolev Institute of Mathematics\\
Celebrates its Fiftieth Anniversary
}
\author{Victor Alexandrov}
\date{Februry 26, 2007}
\maketitle


On May 18, 1957, the Soviet government approved the initiative of
Academicians M.A. Lavrent'ev (1900--1980), S.L. Sobolev
(1908--1989), and S.A. Khristianovich (1908--2000)  aimed at
creating of a new type of research center in Siberia which should
integrate research institutes of all basic scientific,
technological, and humanitarian disciplines, such as mathematics,
physics, mechanics, chemistry, geology, biology, history,
economics, etc.

It was decided to build the center near Novosibirsk, approximately
3000 km east from Moscow. The center has the status of a Branch of
the Academy of Sciences of the USSR and Academician M.A.
Lavrent'ev was appointed as its head.

During 5--10 years,  24 research institutes\,\footnote{During next
50 years the number of institutes was nearly doubled.}, the
Novosibirsk State University, and numerous apartment blocks and
cottages for researchers and staff were built in a picturesque
pine forest on the coast of the manmade lake.

It was the beginning of the famous Akademgorodok (that means
Academy town), to which songs and books are devoted \cite{La82},
which was build by the generation of enthusiasts, devotedly
trusted in the triumph of sciences and human intellect: their
fathers were victors over fascism, their brothers launched the
first sputnik and the first astronaut.

Within the framework of that ambitious project, the Institute of
Mathematics was opened in 1957. The founding father and the first
director of the Institute was Academician Sergej
Sobolev\,\footnote{In 1986 Academician M.M. Lavrent'ev (a son of
M.A. Lavrent'ev) was named his successor, followed by Academician
Yu.L. Ershov who is on duty since 2002.}, one of the most
prominent mathematicians of the 20th century \cite{Le90}.

The main idea was to invite prominent mathematicians from Moscow
and Leningrad\,\footnote{Now St. Petersburg.}, who were willing to
move to Siberia together with their disciples. This idea was
successfully realized. Let us list just a few members of the
Academy of Sciences of the USSR\,\footnote{Now Russian Academy of
Sciences.}, who have worked within the Institute's walls for years
or decades:

$\bullet$ A.D. Alexandrov (1912--1999): the greatest Russian
geometer of the 20th century, the founder of the Soviet school of
geometry `in the large'; who is known world-wide due to his
contributions to the theory of mixed volumes and the theory of
surfaces `in the large', the theory of metric spaces with bounded
curvature and the theory of Monge--Amp\`ere equations, the maximum
principle for elliptic partial differential equations and the
foundations of relativity \cite{Al96}.

$\bullet$ L.V. Kantorovich (1912--1986): a Nobel Prize winner in
economics (1975), one of the creators of a mathematical approach
to economics based on the study of linear extremal problems; his
investigations in the functional analysis, computational
mathematics, the theory of extremal problems, the descriptive
theory of functions and set theory strongly affected those
subjects and gave rise to new fields of research \cite{Ka96}.

$\bullet$ A.A. Lyapunov (1911--1973): starting with descriptive
set theory under the supervision of N.N. Luzin (1883--1950), later
he worked on mathematical aspects of cybernetics and linguistics;
he was awarded the medal `Computer Pioneer' from IEEE Computer
Society (1996).

$\bullet$ A.I. Mal'tsev (1909--1967): the founder of the Siberian
school of algebra and logic; his contributions were mainly  to
algebra (group theory, theory of rings, topological algebra),
mathematical logic (theory of algorithms) and its applications to
algebra \cite{Ma76}.

$\bullet$ S.L. Sobolev (1908--1989): he contributed mainly to the
theory of waves in solids, the theory of equations of mathematical
physics, the functional analysis, the theory of cubature formulas;
he introduced a new class of functional spaces, which are now
known as Sobolev spaces, and the notion of a generalized solution
to a partial differential equation \cite{So06}.

In the early 1990's, the Institute of Mathematics was named after
S.L. Sobolev and, since that time, has been called the Sobolev
Institute of Mathematics or SIM, for short.

In the beginning of 2007 there were 282 research fellows at SIM,
among them 9 members of the Russian Academy of Sciences, 108
professors and 165 fellows with PhD degrees\,\footnote{In fact SIM
has its department in Omsk (the next city with more that 1 million
population west from Novosibirsk, 700 km apart), which
additionally includes 9 professors and 27 fellows with PhD
degrees.}.

\begin{figure}
\begin{center}
\includegraphics[width=4.7 cm]{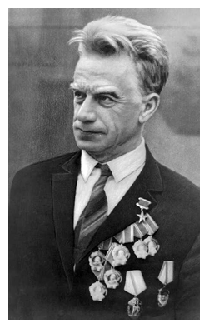}\\
{\bf Sergej Sobolev}
\end{center}
\end{figure}

In general, SIM fellows devote themselves to pure research,
without having obligations to spend time on undergraduate teaching
though many of them supervise post-graduate students and, as a
part-time job, give lectures or even teach undergraduate students
at the Novosibirsk State University.

SIM fellows work in most of the fields of modern mathematics. In
order to provide the reader with an impression of how wide the
variety of research is, we just list a few groups headed by
members of the Russian Academy of Sciences and mention some of
their latest books:

$\bullet$ mathematical logic (Yu.L. Ershov \cite{Er01} and S.S.
Goncharov \cite{Go97}),

$\bullet$ group theory (V.D. Mazurov \cite{Ma06}),

$\bullet$ real functions, potential theory, geometry (Yu.G.
Reshetnyak \cite{Re94}),

$\bullet$ partial differential equations (M.M. Lavrent'ev
\cite{La86} and V.G. Romanov \cite{Ro05}),

$\bullet$ dynamical systems (I.A. Tajmanov \cite{Ta02}),

$\bullet$ probability theory and statistics (A.A. Borovkov
\cite{Bo98}),

$\bullet$ numerical analysis (S.K. Godunov \cite{Go98}).

Approximately 25 research seminars work permanently at SIM. Every
year from 2 to 4 international conferences are organized.

\begin{figure}
\begin{center}
\includegraphics[width=10 cm]{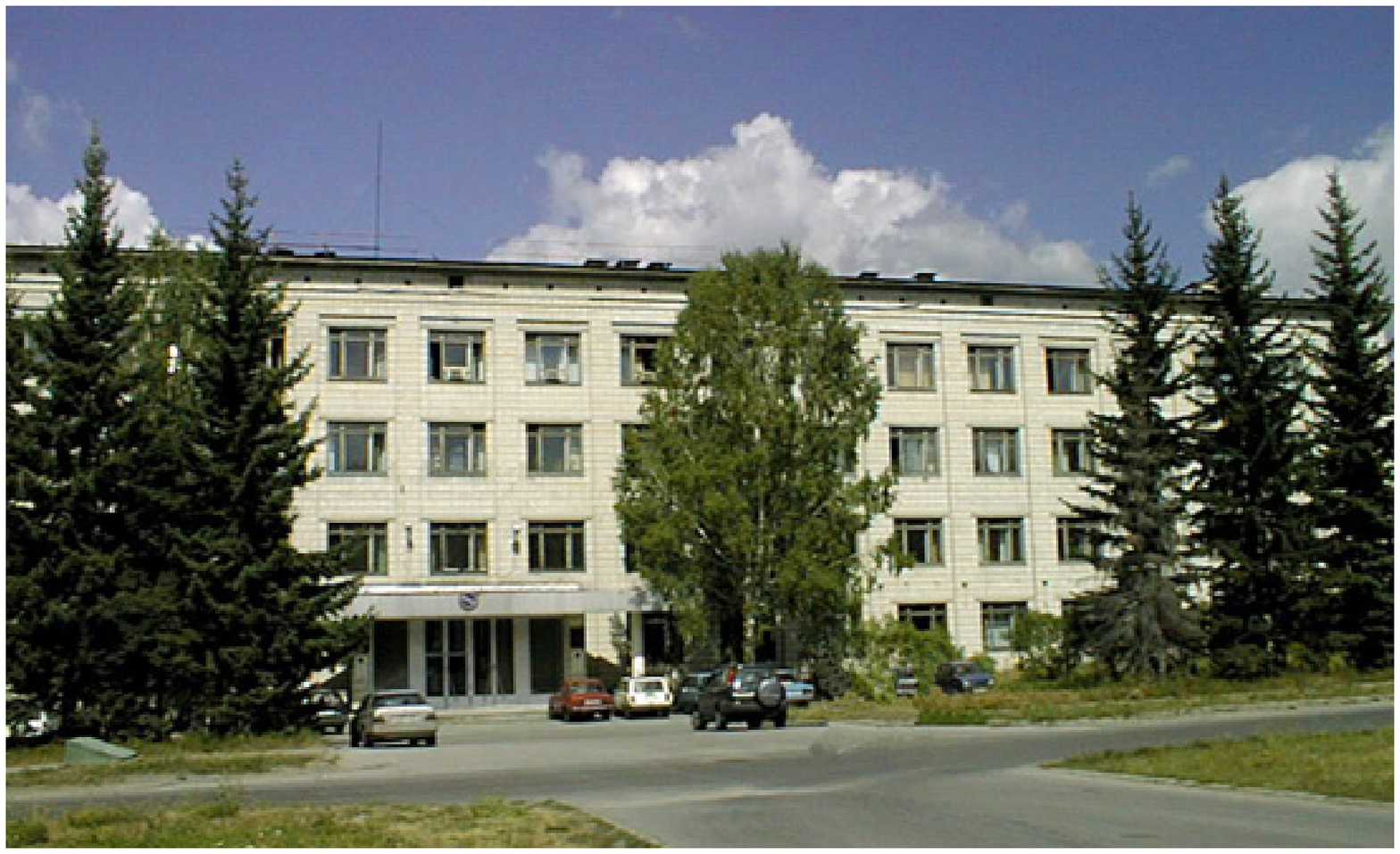}\\
{\bf Sobolev Institute of Mathematics}
\end{center}
\end{figure}

SIM has the right to award the PhD degree or habilitation in the
following fields: mathematical logic, algebra and number theory;
mathematical analysis; geometry and topology; differential
equations; computational mathematics. Each year SIM enrols 12 post
graduate students\,\footnote{The Omsk department enrolls 4
students in addition.} in these fields. These students are
supposed to complete their PhD theses in 3 years.

The SIM library is one of the best mathematical libraries in
Russia East of the Ural mountains. It contains approximately
150,000 items: more that 30,000 books (including  about 20,000
books in foreign languages, including a few books published in the
17th century) and more than 100,000 issues of journals (including
about 75,000 issues of foreign journals).

SIM publishes several journals on mathematics and applied
mathematics in Russian, {\it Algebra and Logic}\,\footnote{ISSN
0373--9252, a cover-to-cover English translation is available.},
{\it Discrete Analysis and Operations Research}\,\footnote{ISSN
1560--7542 for Series I and ISSN 1560--9901 for Series II.}, {\it
Mathematical Transactions}\,\footnote{ISSN 1560--750X, for most of
the articles an English translation is available in {\it Siberian
Advances in Mathematics}, ISSN 1055--1344.}, {\it The Siberian
Journal for Industrial Mathematics}\,\footnote{ISSN 1560--7518.},
{\it Siberian Mathematical Journal}\,\footnote{ISSN 0037--4474, a
cover-to-cover English translation is available.}, and {\it
Siberian Electronic Mathematical Reports}\,\footnote{Electronic
only, available at {\tt http://semr.math.nsc.ru/english.html.}}.

SIM is involved in numerous Russian and international research
programs; it has been a partner of Zentralblatt MATH for more then
10 years. Many mathematicians, who started their careers as SIM
fellows, received international honors for their contributions to
mathematics and received professorships all over the world. To
name a few of them we mention Efim Zel'manov (was awarded a Fields
Medal in 1994; now he is a professor at the University of
California, San Diego); Ivan Shestakov (awarded the 2007 Moore
Research Article Prize; now he is a professor at the University of
Sao Paulo, Brazil); Mikhail Batanin (Macquarie University,
Australia); Oleg Bogopol'skij (Universit\"at Dortmund, Germany);
Serguei Foss (Heriot-Watt University, U.K.); Alexander Kostochka
(University of Illinois, Urbana); Igor Nikolaev (University of
Illinois at Urbana-Champaign);  Vladimir Vershinin (Universit\'e
Montpellier, France); Andrei Voronkov (University of Manchester,
U.K.).

Detailed information can be found on the Institute's web site {\tt
http://math.nsc.ru/english.html.}

\vfill \pagebreak


\begin{thebibliography}{99}

\bibitem{La82}
{\sl  Lavrent'ev, M.A.}
... Will be enlarged by Siberia (in Russian).
Novosibirsk, 1982.

\bibitem{Le90}
{\sl  Leray, J.}
La vie et l'oeuvre de Serge Sobolev.
C. R. Acad. Sci., Paris, S\'er. G\'en., Vie Sci. 7, No.6, 467--471 (1990).

\bibitem{Al96}
{\sl  Alexandrov, A.D.} Selected works. Part 1: Amsterdam: Gordon
and Breach Publishers, 1996. Part 2: Boca Raton, FL: Chapman \&
Hall/CRC, 2005.

\bibitem{Ka96}
{\sl Kantorovich, L.V.} Selected works. Parts 1, 2. Amsterdam:
Gordon and Breach Publishers, 1996.

\bibitem{Ma76}
{\sl Mal'tsev, A.I.} Selected works (in Russian). Vol. 1, 2.
Moscow: Nauka, 1976.

\bibitem{So06}
Selected works of S.L. Sobolev. New York: Springer, 2006.

\bibitem{Er01}
{\sl Ershov, Yu.L.} Multi-valued fields. New York: Kluwer
Academic,  2001.

\bibitem{Go97}
{\sl Goncharov, S.S.} Countable Boolean algebras and decidability.
New York: Plenum, 1997.

\bibitem{Ma06}
{\sl Mazurov, V.D.} (ed.) and {\sl Khukhro, E.I.} (ed.) The
Kourovka notebook. Unsolved problems in group theory. 16th ed.
Novosibirsk: Sobolev Institute of Mathematics, 2006.

\bibitem{Re94}
{\sl Reshetnyak, Yu.G.} Stability theorems in geometry and
analysis. Dordrecht: Kluwer Academic Publishers, 1994.

\bibitem{La86}
{\sl Lavrent'ev, M.M., Romanov, V.G.}, and {\sl Shishatskij, S.P.}
Ill-posed problems of mathematical physics and analysis.
Providence: American Mathematical Society, 1986.

\bibitem{Ro05}
{\sl Romanov, V.G.} Stability in inverse problems (in Russian).
Moscow: Nauchnyj Mir, 2005.

\bibitem{Ta02}
{\sl Tajmanov, I.A.} Lectures in differential geometry (in
Russian). Izhevsk: Institut Komp'yuternykh Issledovanij, 2002.

\bibitem{Bo98}
{\sl Borovkov, A.A.} Probability theory. Abingdon: Gordon and
Breach, 1998.

\bibitem{Go98}
{\sl Godunov, S.K.} Modern aspects of linear algebra. Providence:
American Mathematical Society, 1998.

\end{thebibliography}
\end{document}